\documentclass[12pt,a4paper]{amsart}

\usepackage{amscd, amssymb}

\begin{document}

\pagestyle{headings}

\font\fiverm=cmr5 
\input{prepictex}
\input{pictex}
\input{postpictex}

\def\clsp{\overline{\operatorname{span}}}
\def\newspan{{\operatorname{span}}}

\newtheorem{thm}{Theorem}[section]
\newtheorem{cor}[thm]{Corollary}
\newtheorem{lem}[thm]{Lemma}
\newtheorem{prop}[thm]{Proposition}
\newtheorem{thm1}{Theorem}

\theoremstyle{definition}
\newtheorem{dfn}[thm]{Definition}
\newtheorem{dfns}[thm]{Definitions}

\theoremstyle{remark}
\newtheorem{rmk}[thm]{Remark}
\newtheorem{rmks}[thm]{Remarks}
\newtheorem{example}[thm]{Example}
\newtheorem{examples}[thm]{Examples}
\newtheorem{note}[thm]{Note}
\newtheorem{notes}[thm]{Notes}

\title{$C^*$-algebras of Labelled graphs}
\author{Teresa Bates \& David Pask}
\date{\today}

\begin{abstract}
We describe a class of $C^*$-algebras which simultaneously
generalise the ultragraph algebras of Tomforde and the shift space
$C^*$-algebras of Matsumoto. In doing so we shed some new light on
the different $C^*$-algebras that may be associated to a shift
space. Finally, we show how to associate a simple $C^*$-algebra to
an irreducible sofic shift.

\vspace{3mm} \noindent{\sc keywords:} $C^*$-algebras, labelled
graph, ultragraph, Matsumoto algebra, shift space.

\vspace{3mm} \noindent {\sc MSC (2000):} 46L05, 37B10.
\end{abstract}

\address{School of Mathematics, UNSW, Sydney 2052,
Australia} \email{teresa@unsw.edu.au}

\address{Mathematics, University of Newcastle, NSW 2308, Australia}
\email{david.pask@newcastle.edu.au}

\date{March 17, 2005}
\thanks{This research was supported by the University of Newcastle and the University of New South Wales}

\maketitle

\section{Introduction}

\noindent The purpose of this paper is to introduce a class of
$C^*$-algebras associated to labelled graphs. Our motivation is to
provide a common framework for working with the ultragraph
algebras of Tomforde (see \cite{t,t2}) and the $C^*$-algebras
associated to shift spaces studied by Matsumoto and Carlsen (see
\cite{m,m2,c,cm} amongst others). Here a labelled graph $(E , \pi
)$ over an alphabet ${\mathcal A}$  is a directed graph $E$,
together with a map $\pi : E^1 \to \mathcal{A}$. An ultragraph
$\mathcal{G}$  is a particular example of a labelled graph (see
Example \ref{lgex} (ii)), and a shift space $\Lambda$ has many
presentations as a labelled graph (see \cite{lm}, Example
\ref{lgex} (iii)). Hence it is natural to give our common
framework in terms of labelled graphs.

To a two-sided shift space $\Lambda$ over a finite alphabet,
Matsumoto associates two $C^*$-algebras ${\mathcal O}_\Lambda$ and
${\mathcal O}_{\Lambda^*}$ generated by partial isometries (see
\cite{cm}). Although ${\mathcal O}_\Lambda$ and ${\mathcal
O}_{\Lambda^*}$ are generated by elements satisfying the same
relations, it turns out that they are not isomorphic in general
(see \cite[Theorem 4.1]{cm}). This fact manifests itself in our
realisation in section \ref{ma} of ${\mathcal O}_\Lambda$ and
${\mathcal O}_{\Lambda^*}$ as the $C^*$-algebras of the labelled
graphs $(E_\Lambda,\pi_\Lambda)$ and
$(E_{\Lambda^*},\pi_{\Lambda^*})$ respectively, which are not
necessarily isomorphic as labelled graphs. Moreover, in Corollary
\ref{soficsimple} we show that using labelled graphs gives us the
facility to canonically associate a simple $C^*$-algebra to an
irreducible sofic shift (cf. \cite{cm, c, c1}).

In fact we can associate a number of (possibly different)
$C^*$-algebras to a labelled graph.  This leads us to the notion
of a labelled space, which we describe in section \ref{LS}.
Briefly, a labelled space $( E , \pi , \mathcal{C} )$ consists of
a labelled graph $(E , \pi )$ together with a collection
$\mathcal{C} \subseteq 2^{E^0}$ which plays the same role as
$\mathcal{G}^0$ in \cite{t} and is related to the abelian
AF-subalgebra $A_\Lambda$ (resp. $A_{\Lambda^*}$) in ${\mathcal
O}_\Lambda$ (resp. ${\mathcal O}_{\Lambda^*}$) generated by the
source projections.

In section \ref{CoLS} we define a representation of a labelled space in terms of
partial isometries $\{ s_a : a \in {\mathcal A} \}$ and projections
$\{ p_A : A \in \mathcal{C} \}$ subject to certain relations. Our relations
generalise those found in \cite{t,m}. In order to build a nondegenerate
$C^*$-algebra from a representation of $( E , \pi , \mathcal{C} )$ it is necessary
 for $\mathcal{C}$ to be weakly left-resolving: a condition which is a
generalisation of the left-resolving property for labelled graphs.
Hence we may define $C^* ( E , \pi , \mathcal{C} )$ to be the
$C^*$-algebra which is universal for representations of the weakly
left-resolving labelled space $(E , \pi , \mathcal{C} )$.  Since
any ultragraph has a natural realisation as a left-resolving
labelled graph, the class of $C^*$-algebras of labelled spaces
contains the ultragraph algebras (and hence, graph algebras and
Exel-Laca algebras).

In section \ref{GIUT} we give a version of the gauge invariant uniqueness theorem
for $C^* (E , \pi , \mathcal{C} )$ which will ultimately allow us to make the
connection with the Matsumoto algebras.

In section \ref{App} we give three applications of our uniqueness
theorem: In section \ref{dualgraphs} we show how to construct a
dual labelled space, which is the analogue of the higher block
presentation of a shift space (cf.\ \cite{lm}). We give an
isomorphism theorem for dual labelled spaces  which is a
generalisation of \cite[Corollary 2.5]{bprsz} and forms a starting
point for future work (see \cite{bp3}). In section \ref{ma} we
show that if $\mathcal{O}_\Lambda$ (resp.\
$\mathcal{O}_{\Lambda^*}$) has a gauge action, then it is
isomorphic to the $C^*$-algebra of a certain labelled space. Then
in section \ref{finitenessconditions} we give necessary conditions
for the $C^*$-algebra of a labelled space to be isomorphic to the
$C^*$-algebra of the underlying directed graph. We then show how
to associate a simple $C^*$-algebra to an irreducible shift space.
By example, we show that in general the $C^*$-algebra of a
labelled space will not be isomorphic to the $C^*$-algebra of any
directed graph; hence labelled graph $C^*$-algebras form a
strictly larger class of $C^*$-algebras than graph algebras.

Since we seek to generalise them, we begin by giving a brief description of
Ultragraph algebras and Matsumoto algebras.

\section{Ultragraph Algebras and Matsumoto Algebras} \label{niceguys}

\subsection{Ultragraph Algebras}

\noindent An ultragraph $\mathcal{G} = ( G^0  , \mathcal{G}^1 , r
, s )$ consists of a countable set of vertices $G^0$, a countable
set of edges $\mathcal{G}^1$, and functions $s : \mathcal{G}^1 \to
G^0$ and $r : \mathcal{G}^1 \to 2^{G^0}$. Let $\mathcal{G}^0$ be
the smallest collection of $2^{G^0}$ which contains $s(e)$ and
$r(e)$ for all $e \in \mathcal{G}^1$ and is closed under finite
intersections and unions. The ultragraph algebra $C^* (
\mathcal{G} )$ is the universal $C^*$-algebra for Cuntz-Krieger
$\mathcal{G}$-families: collections of partial isometries $\{ s_e
: e \in \mathcal{G}^1 \}$ with mutually orthogonal ranges, and
projections $\{ p_A : A \in \mathcal{G}^0 \}$ satisfying the
relations
\begin{itemize}
\item[1.]  $p_\emptyset = 0$, $p_A p_B = p_{A \cap B}$ and $p_{A
\cup B} = p_A + p_B - p_{A \cap B}$ for all $A,B \in {\mathcal
G}^0$

\item[2.] $s_e^*s_e = p_{r(e)}$ and $s_e s_e^* \leq p_{s(e)}$ for
all $e \in {\mathcal G}^1$

\item[3.]  $p_v = \sum_{s(e) = v} s_e s_e^*$ whenever $0 <
|s^{-1}(v)| < \infty$
\end{itemize}

\noindent (see \cite[Definition 2.7]{t}). Recall that $v \in G^0$
is an infinite emitter if $\vert s^{-1} (v) \vert = \infty$. We
say that $\mathcal{G}$ is {\em row-finite} if there are no
infinite emitters and $r(e)$ is finite for all $e \in
\mathcal{G}^1$.

Ultragraph algebras simultaneously generalise graph $C^*$-algebras
and Exel-Laca algebras (see \cite[Sections 3 and 4]{t}). By
\cite[Corollary 5.5]{t2} there is a non row-finite ultragraph
whose $C^*$-algebra is not Morita equivalent to a graph algebra or
an Exel-Laca algebra.

\subsection{Matsumoto Algebras}

For an introduction to shift spaces we refer the reader to the
excellent treatment in \cite{lm}. Let $\Lambda$ be a two-sided
shift space over a finite alphabet ${\mathcal A}$ and denote by
$\Lambda^*$ the collection of all finite blocks in $\Lambda$. Let
\begin{equation} \label{xlambdadef}
X_\Lambda = \{ ( x_i )_{i \ge 1} : ( x_i )_{i \in {\bf Z}} \in
\Lambda \}
\end{equation}

\noindent  denote the set of all right-infinite sequences in
$\Lambda$.

Following \cite{cm} there are two $C^*$-algebras associated to
$\Lambda$. Each $C^*$-algebra is
generated by partial isometries $\{ t_a : a \in \mathcal{A} \}$ subject to
\begin{equation} \label{mrels}
\sum_{a \in \mathcal{A}} t_a t_a^* = 1 , \text{ and }
t_\alpha^* t_\alpha t_\beta = t_\beta t_{\alpha \beta}^* t_{\alpha \beta} ,
\text{ where } \alpha , \beta , \alpha \beta \in \Lambda^* .
\end{equation}

\noindent As in \cite{cm} we denote by $\mathcal{O}_\Lambda$ the
$C^*$-algebra defined directly on Hilbert space in
\cite{m4,m45,m5} and by $\mathcal{O}_{\Lambda^*}$ the
$C^*$-algebra defined using the Fock space construction in
\cite{m,m2,m1,m3,m6}. Because of the different ways in which the
relations (\ref{mrels}) are realised it turns out that
$\mathcal{O}_\Lambda$ and $\mathcal{O}_{\Lambda^*}$ are not
isomorphic in general (see \cite[Section 6]{cm}).


There is a uniqueness theorem for ${\mathcal O}_\Lambda$ (resp.\
$\mathcal{O}_{\Lambda^*}$) when $\Lambda$ satisfies condition (I)
(resp.\ condition $(I^*)$)  given in \cite[Section 4]{cm} (resp.\
\cite[Section 3]{cm}).

{\bf Condition $(I)$:} For $x \in X_\Lambda$ and $l \in {\bf N}$
put $\Lambda_l (x) = \{ \mu \in \Lambda_l : \mu x  \in X_\Lambda
\}$. Two infinite paths $x , y \in X_\Lambda$ are {\em $l$-past
equivalent} (written $x \thicksim_l y$) if $\Lambda_l (x) =
\Lambda_l (y)$. The shift space $X_\Lambda$ satisfies condition
(I) if for any $l \in {\bf N}$ and $x \in X_\Lambda$ there exists
$y \in X_\Lambda$ such that $y \neq x$, $y \sim_l x$.

{\bf Condition $(I^*)$:} For $\omega \in \Lambda^*$ and $l \in
{\bf N}$ we set $\Lambda_l(\omega) = \{ \mu :  |\mu| \leq l, \;
\mu \omega \in \Lambda^*\}$. Two words $\mu, \nu \in \Lambda^*$
are said to be {\em $l$-past equivalent} (written $\mu \thicksim_l
\nu$ ) if $\Lambda_l(\mu) = \Lambda_l(\nu)$. The subset
$\Lambda^*_l \subseteq \Lambda^*$ is defined by
$$
\Lambda^*_l := \{ \omega \in \Lambda^* : |\{\mu \in \Lambda^* : \mu \thicksim_l
\omega \}| <  \infty \}.
$$
\noindent
The shift space $\Lambda$ satisfies condition $(I^*)$ if for every
$l \in {\bf N}$ and $\mu \in \Lambda_l^*$ there exist
distinct words $\xi_1,\xi_2 \in \Lambda^*$ with $|\xi_1| = |\xi_2| = m$
such that
$$
\mu \thicksim_l \xi_1 \gamma_1 \mbox{  and  } \mu \thicksim_l \xi_2 \gamma_2
$$
for some $\gamma_1,\gamma_2 \in \Lambda^*_{l + m}$.

\begin{thm} \label{matsumotohasgauge}
Let $\Lambda$ be a two-sided shift space over a finite alphabet
which satisfies condition (I) (resp.\ $(I^*)$).  Then there is a
strongly continuous action $\beta$ (resp.\ $\beta^*$) of ${\bf T}$
on ${\mathcal O}_\Lambda$ (resp.\ $\mathcal{O}_{\Lambda^*}$) such
that $\beta_z (t_a) = z t_a$ (resp.\ $\beta^*_z (t_a ) = z t_a$)
for all $a \in \mathcal{A}$ and $z \in {\bf  T}$.
\end{thm}

\begin{proof}
That each $\beta_z$ (resp.\ $\beta^*_z$) is an automorphism of
$\mathcal{O}_\Lambda$ (resp.\ $\mathcal{O}_{\Lambda^*}$) for each
$z \in {\bf T}$ follows from
\cite[Proposition 4.2]{cm} (resp.\ \cite[Proposition 3.4]{cm}).
A standard $\frac{\epsilon}{3}$ argument shows that $\beta$ (resp.\ $\beta^*$)
is strongly continuous.
\end{proof}

\noindent In \cite{m6} Matsumoto defines $\lambda$-graph systems
${\mathcal L}_\Lambda$ and ${\mathcal L}_{\Lambda^*}$ associated
to a two-sided shift space $\Lambda$ together with corresponding
$C^*$-algebras ${\mathcal O}_{{\mathcal L}_\Lambda}$ and
 ${\mathcal O}_{{\mathcal L}_{\Lambda^*}}$.  By \cite[Theorem
5.6]{cm} we see that if $\Lambda$ satisfies condition (I) then
${\mathcal O}_\Lambda \cong  {\mathcal O}_{{\mathcal L}_\Lambda}$
and if $\Lambda$ satisfies condition $(I^*)$ then ${\mathcal
O}_{\Lambda^*} \cong  {\mathcal O}_{{\mathcal L}_{\Lambda^*}}$.
Hence, for our purposes, it suffices to work with ${\mathcal
O}_\Lambda$ and ${\mathcal O}_{\Lambda^*}$.

\section{Labelled spaces} \label{LS}

\noindent A directed graph $E$ consists of a quadruple $( E^0 , E^1
, r , s )$ where $E^0$ and $E^1$ are countable sets of vertices
and edges respectively and $r, s : E^1 \to E^0$ are maps giving
the direction of each edge. A path $\lambda = e_1 \ldots e_n$ is a
sequence of edges $e_i \in E^1$ such that $r ( e_i ) = s ( e_{i+1}
)$ for $i=1 , \ldots , n-1$. The collection of paths of length $n$
in $E$ is denoted $E^n$ and the collection of all finite paths in
$E$ by $E^*$, so that $E^* = \cup_{n \ge 0} E^n$.
The edge shift $( \textsf{X}_E , \sigma_E )$ associated
to a directed graph  $E$  with no sinks or sources is defined by:
\[
\textsf{X}_E = \{ x \in ( E^1 )^{\bf Z} : s ( x_{i+1} ) = r (
x_{i} ) \text{ for all } i \in {\bf Z} \} \text{ and } ( \sigma_E
x )_i = x_{i+1} \text{ for } i \in {\bf Z} .
\]

\noindent
The following definition is adapted from \cite[Definition 3.1.1]{lm}:

\begin{dfn}
A {\it labelled graph} $( E , \pi )$ over an alphabet ${\mathcal
A}$ consists of a directed graph $E$ together with a labelling map
$\pi : E^1 \to \mathcal{A}$.
\end{dfn}

\noindent Without loss of generality we may assume that
the map $\pi$ is onto. Given a labelled graph $(E, \pi)$ such that
every vertex in $E$ emits and receives an edge, we may define a subshift
$( \textsf{X}_{(E,\pi)} , \sigma )$ of $\mathcal{A}^{\bf Z}$ by
\[
\textsf{X}_{(E,\pi)} = \{ y \in \mathcal{A}^{\bf Z} : \text{ there exists
 } x \in \textsf{X}_E \mbox{ such that } y_i = \pi ( x_i ) \text{ for all } i \in {\bf Z}
\} ,
\]

\noindent where $\sigma$ is the shift map. The labelled graph $(E
, \pi )$ is said to be a {\em presentation} of the shift space $X
= \textsf{X}_{(E , \pi )}$. As shown in \cite[\S 3.1]{lm} a shift
space may have many different presentations (see Examples
\ref{lgex} (ii), (vi), (vii)).

Let $\mathcal{A}^*$ be the collection of all {\em words} in the
symbols of $\mathcal{A}$ (see \cite[\S 0.2]{Ta}). The map $\pi$
extends naturally to a map $\pi : E^n \to \mathcal{A}^*$, where $n
\ge 1$: for $\lambda = e_1 \ldots e_n \in E^n$ put $\pi ( \lambda
)  = \pi ( e_1 ) \ldots \pi ( e_n )$; in this case the path
$\lambda \in E^n$ is said to be a {\em representative} of the {\em
labelled path} $\pi ( e_1 ) \ldots \pi ( e_n )$. Let
$\mathcal{L}^n ( E , \pi ) = \pi ( E^n )$ denote the collection of
all labelled paths in $(E,\pi)$ of length $n$, then $\mathcal{L}
(E,\pi) = \cup_{n \ge 1} \mathcal{L}^n (E,\pi)$ denotes the
collection of all words in the alphabet $\mathcal{A}$ which may be
represented by paths in the labelled graph $(E, \pi )$. In this
way $\pi$ induces a map from the language $\cup_{n \ge 1} E^n$ of
the subshift of finite type $\textsf{X}_E$ associated to $E$ into
$\mathcal{L} (E,\pi)$, the language of the shift space
$\textsf{X}_{(E,\pi)}$ presented by $(E,\pi)$ (see \cite[\S
3]{lm}). The usual length function $\vert \cdot \vert : E^* \to
\bf{N}$ transfers naturally over to $\mathcal{L} ( E, \pi)$.

For $\alpha$ in $\mathcal{L} ( E , \pi )$ we put
\[
s_\pi ( \alpha ) = \{ s ( \lambda ) \in E^0 : \pi ( \lambda ) =
\alpha \} \text{ and } r_\pi ( \alpha ) = \{ r ( \lambda ) \in E^0
: \pi ( \lambda ) = \alpha \} ,
\]

\noindent so that $r_\pi , s_\pi : \mathcal{L} ( E , \pi ) \to
2^{E^0}$. We shall drop the subscript on $r_\pi$ and $s_\pi$ if
the context in which it is being used is clear. For $\alpha ,
\beta \in \mathcal{L} (E,\pi)$ we have $\alpha \beta \in
\mathcal{L} ( E , \pi )$ if and only if $r ( \alpha ) \cap s (
\beta ) \neq \emptyset$.

Where possible we shall denote the elements of $\mathcal{A} =
\mathcal{L}^1 (E , \pi)$ as $a,b$, etc., elements of $\mathcal{L}
(E , \pi )$ as $\alpha , \beta$, etc., leaving $e, f$ for elements
of $E^1$ and $\lambda  , \mu$ for elements of $E^*$.

Let $(E , \pi )$ and $(F , \pi' )$ be graphs labelled by the same
alphabet. A graph isomorphism $\phi : E \to F$ is a
{\em labelled graph isomorphism} if $\pi' ( \phi (e) ) = \pi (e)$ for
all $e \in E^1$ and we write $\phi : ( E , \pi ) \to ( F , \pi'
)$.

\begin{dfn}
The labelled graph $( E , \pi )$ is {\em left-resolving} if for
all $v \in E^0$ the map $\pi : r^{-1} (v) \to \mathcal{A}$ is
injective.
\end{dfn}

\noindent The left-resolving condition ensures that for all $v \in
E^0$ the labels $\{ \pi (e) : r(e) = v \}$ of all incoming edges
to $v$ are all different. When $( E , \pi )$ is left-resolving, if
$\lambda , \mu \in \cup_{n \ge 1} E^n$ satisfy $\pi ( \lambda ) =
\pi ( \mu )$ and $r ( \lambda ) = r ( \mu )$ then $\lambda = \mu$.

\begin{examples} \label{lgex}
\begin{itemize}

\item[(i)] Let $E$ be a directed graph. Put $\mathcal{A} = E^1$
and let $\pi : E^1 \to E^1$ be the identity map (the
\textit{trivial labelling}) then $( E , \pi )$ is a left-resolving
labelled graph.

\item[(ii)] Let $\mathcal{G} = ( G^0 , \mathcal{G}^1 , r , s )$ be
an ultragraph. Define $E = E_{\mathcal{G}}$ by putting $E^0 =
G^0$, $E^1 = \{ ( s(e) , w ) : e \in \mathcal{G}^1 , w \in r (e)
\}$ and defining $r' , s' : E^1 \to E^0$ by $s' ( v , w ) = v$,
$r' ( v , w ) = w$. Set $\mathcal{A}=\mathcal{G}^1$ and define
$\pi_{\mathcal G} : E^1 \to \mathcal{A}$ by $\pi (s(e),w) = e$.
The resulting labelled graph $(E_{\mathcal{G}},\pi_{\mathcal G})$
is left-resolving since the source map is single-valued. If
$\mathcal{G}$ is row-finite then $E_\mathcal{G}$ is row-finite.

Conversely, given a left-resolving labelled graph $(E , \pi )$
over an alphabet ${\mathcal A}$ where $s_\pi : \mathcal{L} (E ,
\pi ) \to 2^{E^0}$ is single-valued, we can form a ultragraph
$\mathcal{G}_E = ( E^0 , \mathcal{A} , r', s')$ with $s'=s_\pi$
and $r'=r_\pi$. If $E$ is row finite then the ultragraph
$\mathcal{G}_E$ is row-finite.

\item[(iii)] Following \cite[\S 3]{lm} the labelled graphs $
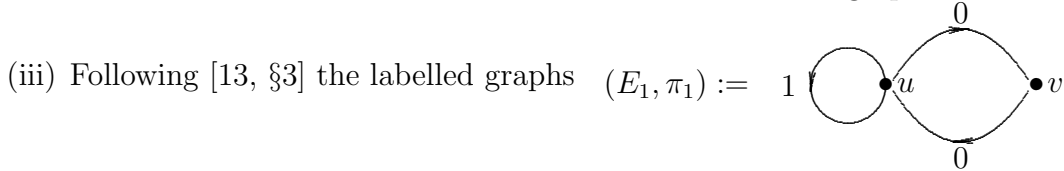
\beginpicture
\setcoordinatesystem units <1cm,0.75cm>

\setplotarea x from -1 to 5, y from -0.2 to 1.1




\put{$(E_1 , \pi_1 ) :=$}[l] at -0.75 0

\put{$\bullet$} at 3 0

\put{$\bullet$} at 5 0

\put{$1$}[l] at 1.6 0

\put{$0$}[b] at 4 1.1

\put{$0$}[t] at 4 -1.1

\put{$u$}[l] at 3.15 0

\put{$v$}[l] at 5.15 0


\setquadratic

\plot 3.1 0.1  4 1 4.9 0.1 /

\plot 3.1 -0.1 4 -1  4.9 -0.1 /

\circulararc 360 degrees from 3 0 center at 2.5 0

\arrow <0.25cm> [0.1,0.3] from 2.015 0.1 to 2 -0.1 

\arrow <0.25cm> [0.1,0.3] from 4.1 -0.985 to 3.9 -1

\arrow <0.25cm> [0.1,0.3] from 3.9 0.985 to 4.1 1

\endpicture
$
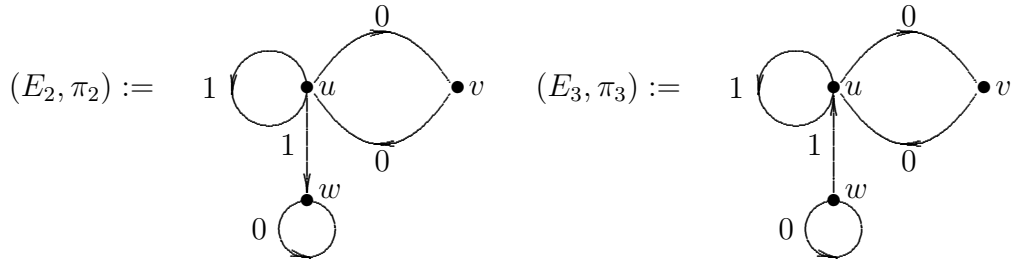
\[
\beginpicture
\setcoordinatesystem units <1cm,0.75cm>

\setplotarea x from -2.5 to 13, y from -3.1 to 1.1




\put{$(E_2 , \pi_2 ) :=$} at 0 0

\put{$\bullet$} at 3 0

\put{$\bullet$} at 5 0

\put{$\bullet$} at 3 -2

\put{$1$}[l] at 1.6 0

\put{$0$}[b] at 4 1.1

\put{$0$}[t] at 4 -1.1

\put{$u$}[l] at 3.15 0

\put{$v$}[l] at 5.15 0

\put{$w$}[bl] at 3.15 -1.95

\put{$1$}[r] at 2.85 -1

\put{$0$}[r] at 2.47 -2.5

\setquadratic

\plot 3.1 0.1  4 1 4.9 0.1 /

\plot 3.1 -0.1 4 -1  4.9 -0.1 /

\circulararc 360 degrees from 3 0 center at 2.5 0

\circulararc 360 degrees from 3 -2 center at 3 -2.5

\arrow <0.25cm> [0.1,0.3] from 2.015 0.1 to 2 -0.1 

\arrow <0.25cm> [0.1,0.3] from 4.1 -0.985 to 3.9 -1

\arrow <0.25cm> [0.1,0.3] from 3.9 0.985 to 4.1 1

\arrow <0.25cm> [0.1,0.3] from 3 -0.15 to 3 -1.85

\arrow <0.25cm> [0.1,0.3] from 2.95 -2.98 to 3.04 -3.01


\put{$(E_3 , \pi_3 ) :=$} at 7 0

\put{$\bullet$} at 10 0

\put{$\bullet$} at 12 0

\put{$\bullet$} at 10 -2

\put{$1$}[l] at 8.6 0

\put{$0$}[b] at 11 1.1

\put{$0$}[t] at 11 -1.1

\put{$u$}[l] at 10.15 0

\put{$v$}[l] at 12.15 0

\put{$w$}[bl] at 10.15 -1.95

\put{$1$}[r] at 9.85 -1

\put{$0$}[r] at 9.47 -2.5

\setquadratic

\plot 10.1 0.1  11 1 11.9 0.1 /

\plot 10.1 -0.1 11 -1  11.9 -0.1 /

\circulararc 360 degrees from 10 0 center at 9.5 0

\circulararc 360 degrees from 10 -2 center at 10 -2.5

\arrow <0.25cm> [0.1,0.3] from 9.015 0.1 to 9 -0.1 

\arrow <0.25cm> [0.1,0.3] from 11.1 -0.985 to 10.9 -1

\arrow <0.25cm> [0.1,0.3] from 10.9 0.985 to 11.1 1

\arrow <0.25cm> [0.1,0.3] from 10 -1.85 to 10 -0.15

\arrow <0.25cm> [0.1,0.3] from 9.95 -2.98 to 10.04 -3.01

\endpicture
\]

\noindent have the same language as the even shift $Y$ since
between any two $1$'s there must be an even number of $0$'s. Hence
$\textsf{X}_{(E_i , \pi_i)} = Y$ for $i=1,2,3$ by
\cite[Proposition 1.3.4 (3)]{lm}. Only graphs $( E_1 ,\pi_1 )$ and
$( E_2 , \pi_2 )$ are left-resolving.

\item[(iv)] Let $E$ be a directed graph and $\Gamma$ a group which
acts on (the right of) $E$. Define $\pi_q : E^1 \to E^1 / \Gamma$
by $\pi_q ( e ) = q (e)$ where $q : E^1 \to E^1 / \Gamma$ is the
quotient map. If the action of $\Gamma$ is free on $E^1$, then the
resulting labelled graph $( E , \pi_q )$ is left-resolving.
More generally, if $p: F \to E$ is a graph morphism then there is a
labelling $\pi_p : F^1 \to E^1$ given by $\pi_p (f) = p(f)$ for
all $f \in E^1$. If $p$ is a covering map then $\pi_p$ is
left-resolving.

\item[(v)] Recall from \cite[\S 3]{bp}, that an out-splitting of a
directed graph $E$ is formed by a partition $\mathcal{P}$ of
$s^{-1} (v)$ into $m(v) \ge 1$ non-empty subsets for each $v \in
E^0$ (if $s^{-1} (v) = \emptyset$ then $m(v)=0$). Given such a
partition $\mathcal{P}$ one may construct a directed graph $E_s
({\mathcal P})$ where $E_s ( \mathcal{P}^1 ) = \{ e^j : e \in E^1
, 1 \le j \le m ( r (e) ) \} \cup \{ e : m (r(e)) = 0 \}$. Define
$\pi : E_s ( \mathcal{P} )^1 \to E^1$ by $\pi ( e^j ) = e$ for $1
\le j \le m ( r (e) )$ and $\pi (e) = e$ if $m (r(e) )=0$. For an
in-splitting (see \cite[\S 5]{bp}) of $E$ using a partition
$\mathcal{P}$, a similar construction also yields a labelled
graph. However the resulting labelling $\pi$ of the in-split graph
$E_r ( \mathcal{P} )$ will not be left-resolving in general.

\item[(vi)] Let $\Lambda$ be a two-sided shift space over a finite
alphabet $\mathcal{A}$ with $X_\Lambda$ defined as in
(\ref{xlambdadef}). Let $X_\Lambda^- = \{ ( x_i )_{i \le 0} : (
x_i )_{i \in {\bf Z}} \in \Lambda \}$ so that any element $x \in
\Lambda$ may be written as $x = x^- x^+$. For arbitrary $x^+ \in
X_\Lambda$ and $x^- \in X_\Lambda^-$ the bi-infinite sequence $y =
x^- x^+$ may not belong to $\Lambda$. Define the past set of $t
\in X_\Lambda$ as
\[
P_{\infty}( t ) = \{ x^- \in X^-_\Lambda \; : \; x^- t \in \Lambda
\}.
\]

\noindent A shift is {\em sofic} if and only if the number of past
sets is finite \cite{Kr,lm}.

For $s ,t \in X_\Lambda$, we say that $s$ is {\em past equivalent}
to $t$ (denoted $s \sim_{\infty} t$) if $P_{\infty}(s) =
P_{\infty}(t)$. Evidently $X_\Lambda / \sim_{\infty}$ can be
identified with $\Omega_\Lambda$ as described in \cite[\S 2]{m2}.
Define a labelled graph $( E_\Lambda , \pi_\Lambda )$ as follows:
let $E_\Lambda^0 = \{ [v] : v \in X_\Lambda / \sim_\infty \}$,
$E_\Lambda^1 = \{ ( [v] , a, [w] ) : a \in \mathcal{A}, a w
\sim_\infty v \}$ with $s( [v], a, [w]) = [v]$ and $r ( [v] ,a,
[w] ) = [w]$. If $([v] ,a, [w] ) \in E_\Lambda^1$ we put
$\pi_\Lambda ( [v] , a, [w] ) = a$. The resulting left-resolving
labelled graph is usually referred to as the {\em left-Krieger
cover} of $\Lambda$ and the construction is evidently independent
of the choice of representatives (see \cite{Kr}).

If $Y$ is the even shift then $( E_Y , \pi_Y )$ is labelled graph
isomorphic to $( E_2 , \pi_2 )$ in (iii) above. Let $Z$ be shift
over the alphabet $\{ 1 , 2, 3 , 4 \}$ in which the words
\[
\{ 1 2^k 1 , 3 2^k 1 2, 3 2^k 13 , 4 2^k 14 : k \ge 0 \}
\]

\noindent do not occur (see \cite[\S 4]{cm}) then $( E_Z , \pi_Z
)$ has six vertices.

\item[(vii)]  Let $\Lambda$ be a two-sided shift over a finite
alphabet $\mathcal{A}$. We construct a variant of the predecessor
graph $(E_{\Lambda^*},\pi_{\Lambda^*})$ in the following way. For
$\mu \in \Lambda^*$ we define
$$
P(\mu) := \{ \lambda : \lambda \mu \in \Lambda^* \}
$$
and define an equivalence relation by $\mu \thicksim \nu$ if
$P(\mu) = P(\nu)$. A shift is sofic if and only if the number of
predecessor sets is finite \cite{lm}.

Let $\Lambda^*_\infty$ denote those $\mu \in \Lambda^*$ which have
an infinite equivalence class. Since $\mathcal{A}$ is finite
$\Lambda^*_\infty / \thicksim$ can be identified with
$\Omega_{\Lambda^*} = \lim_\leftarrow \Omega_l^*$ as described in
\cite[Section 2]{m2}. We set $E^0_{\Lambda^*} = \Lambda^*_\infty /
\thicksim$, $E^1_{\Lambda^*} = \{ ([\mu],a,[\nu]) : a \in
{\mathcal A}, [\mu] = [a \nu] \}$, $r([\mu],a,[\nu]) = [\nu]$ and
$s([\mu],a,[\nu]) = [\mu]$.  The labelling map is defined by
$\pi_{\Lambda^*} ([\mu],a,[\nu])= a$. The resulting labelled graph
is evidently left-resolving.

If $Y$ is the even shift then $( E_{Y^*} , \pi_{Y^*} )$ is
labelled graph isomorphic to $( E_2 , \pi_2 )$ in (iii) above
(cf.\ \cite{c,m2}). If $Z$ is the sofic shift described in Example
\ref{lgex} (vi) then $( E_{Z^*} , \pi_{Z^*} )$ has seven vertices
and contains $( E_Z , \pi_Z )$ as a subgraph.
\end{itemize}
\end{examples}

\begin{dfn} Let $(E,\pi)$ be a labelled graph. For $A \subseteq E^0$
and $\alpha \in \mathcal{L} (E,\pi)$ the {\em relative range of
$\alpha$ with respect to $A$} is defined to be
\[
r(A,\alpha) = \{ r ( \lambda ) : \lambda \in E^* , \pi ( \lambda )
= \alpha , s ( \lambda ) \in A \} .
\]
\end{dfn}

\begin{rmk} \label{relunionintersect}
For any $A, B \subseteq E^0$ we have
\[
r ( A \cap B , \alpha ) \subseteq r ( A , \alpha ) \cap r ( B ,
\alpha ) \text{ and } r ( A \cup B , \alpha ) = r ( A , \alpha )
\cup r ( B , \alpha ) . \]

\noindent For all $A \subseteq {E}^0$ and $\alpha \in \mathcal{L}
(E,\pi)$ we have $r ( A , \alpha ) = r ( A \cap s ( \alpha ) ,
\alpha )$.
\end{rmk}

\noindent A collection $\mathcal{C} \subseteq 2^{E^0}$ of subsets
of $E^0$ is said to be {\em closed under relative ranges for
$(E,\pi)$} if for all $A \in \mathcal{C}$ and $\alpha \in
\mathcal{L} (E,\pi)$ we have $r (A,\alpha) \in \mathcal{C}$. If
$\mathcal{C}$ is closed under relative ranges for $(E , \pi )$,
contains $r ( \alpha )$ for all $\alpha \in \mathcal{L} (E , \pi
)$ and is also closed under finite intersections and unions, then
we say that $\mathcal{C}$ is {\em accommodating} for $(E, \pi)$.

\begin{dfn} A {\em labelled space} consists of a triple $(E , \pi ,
  \mathcal{C} )$,
where $(E , \pi )$ is a labelled graph and $\mathcal{C}$ is
accommodating for $(E , \pi )$ .
\end{dfn}

\begin{dfn} A labelled space $(E , \pi , \mathcal{C} )$ is {\em weakly left-resolving}
if for every $A , B \in \mathcal{C}$ and every $\alpha \in
\mathcal{L} ( E , \pi )$ we have $r ( A , \alpha ) \cap r ( B ,
\alpha ) = r ( A \cap B , \alpha)$.
\end{dfn}

\noindent In particular $(E , \pi , \mathcal{C} )$ is weakly
left-resolving if no pair of disjoint sets $A , B \in \mathcal{C}$
can emit paths $\lambda , \mu$ respectively with $\pi ( \lambda )
= \pi ( \mu )$ and $r ( \lambda ) = r ( \mu )$.  If $( E , \pi )$
is left-resolving then $( E , \pi , \mathcal{C} )$ is weakly
left-resolving for any $\mathcal{C}$.  Evidently if $(E, \pi,
\mathcal{C})$ is weakly left-resolving, then
$(E,\pi,\mathcal{C}')$ is weakly left-resolving for any
$\mathcal{C}' \subseteq \mathcal{C}$.

\noindent Consider the following subsets of $2^{E^0}$
\begin{align*}
\mathcal{E} &= \{ \{ v \} : v \in E^0 \text{ is a source or a sink } \}
\cup \{ r ( \alpha ) : \alpha \in \mathcal{L} (E,\pi) \} \cup \{ s
( \alpha ) : \alpha \in \mathcal{L} (E,\pi) \} \\
\mathcal{E}_- &=\{ \{ v \} : v \in E^0 \text{ is a sink } \}
\cup \{ r ( \alpha ) : \alpha \in \mathcal{L} (E,\pi) \} .
\end{align*}


\noindent The following definition is analogous to the definition
of ${\mathcal G}^0$ in \cite{t}.
\begin{dfn} \label{E0def}
Let $\mathcal{E}^0$ (resp.\ $\mathcal{E}_-^0$) denote the smallest
subset of $2^{E^0}$ containing $\mathcal{E}$ (resp.\
$\mathcal{E}_-$) which is accommodating for $(E,\pi)$ .
\end{dfn}

\begin{rmk} \label{relunionintersect2}
For all $\alpha \in \mathcal{L} (E , \pi )$ we have $r ( s (
\alpha ) , \alpha ) = r ( \alpha )$; moreover if $( \alpha , \beta
) \in \mathcal{L}^{(2)} ( E , \pi )$ then $r ( r ( \alpha ) ,
\beta) = r ( \alpha\beta)$.
For $\alpha , \beta \in \mathcal{L} (E , \pi )$ with $\alpha \beta
\in \mathcal{L}( E , \pi )$ and $A \subseteq {E}^0$ we have $r ( r
( A , \alpha ) , \beta ) = r ( A , \alpha \beta )$.
\end{rmk}

\noindent For labelled spaces $(E,\pi,{\mathcal E}^0)$ which are
weakly left-resolving Remark \ref{relunionintersect} and Remark
\ref{relunionintersect2} show that to form $\mathcal{E}^0$ it
suffices to form
\[
\mathcal{E} \cup \{ r ( A , \alpha ) : A \in \mathcal{E} ,
\alpha \in \mathcal{L} (E , \pi ) \}
\]

\noindent and then close under finite intersections and unions.
To form $\mathcal{E}^0_-$, by Remark \ref{relunionintersect} it suffices to close $\mathcal{E}_-$ under
finite intersections and unions. Evidently, $\mathcal{E}^0_- \subseteq
\mathcal{E}^0$; the containment can be strict, for instance this
occurs when $E$ has sources.
One can show that $\mathcal{E}^0 =
\mathcal{E}^0_-$ if and only if for every $\alpha \in \mathcal{L}(E,\pi)$,
$s ( \alpha )$ can be written as a finite union of
sets of the form $\cap_{i=1}^n r ( \beta_i )$.
Since $E^0$, $\mathcal{L} (E, \pi )$ and $\mathcal{E}$ are
countable it follows that $\mathcal{E}^0$ and $\mathcal{E}_-^0$ are
countable.

For $A \in 2^{{E}^0}$ and $n \ge 1$ let
\[
L_A^n = \{ \alpha \in \mathcal{L}^n ( E , \pi ) : A \cap s (
\alpha ) \neq \emptyset \}
\]

\noindent denote those labelled paths of length $n$ whose source
intersects $A$ nontrivially.

\section{$C^*$-algebras of labelled spaces} \label{CoLS}


\begin{dfn} \label{lgdef}
Let $( E , \pi , \mathcal{C} )$ be a weakly left-resolving labelled space. A
representation of $( E , \pi , \mathcal{C} )$
consists of projections $\{ p_A : A
\in \mathcal{C} \}$ and partial isometries $\{ s_a : a \in
\mathcal{L}^1 ( E , \pi ) \}$ with the properties that

\begin{itemize}

\item[(i)] If $A, B \in \mathcal{C}$ then $p_A p_B = p_{A \cap
B}$ and $p_{A \cup B} = p_A + p_B - p_{A \cap B}$, where
$p_\emptyset = 0$.

\item[(ii)] If $a \in \mathcal{L}^1 (E,\pi)$ and $A \in
\mathcal{C}$ then $p_A s_a = s_a p_{r ( A, a )}$.

\item[(iii)] If $a , b \in \mathcal{L}^1 (E , \pi )$ then $s_a^*
s_a = p_{r( a )}$ and $s_{a}^*s_{b} = 0$ unless $a = b$.

\item[(iv)] 
For $A \in \mathcal{C}$, if $L^1_A$ is
finite and non-empty we have
\begin{equation} \label{sumcond}
p_A = \sum_{a \in L^1_A
} s_{a} p_{r( A , a )} s_{a}^* .
\end{equation}
\end{itemize}

\end{dfn}

\noindent If $a , b \in \mathcal{L}^1 (E , \pi )$ are such that $a
b \in \mathcal{L} ( E , \pi )$ then we have
\[
( s_a^* s_a )( s_b s_b^* ) = p_{r ( a )} s_b s_b^* = s_b p_{r ( r
( a ) , b )} s_b^* = s_b s_b^* p_{r ( a )} = ( s_b s_b^* )( s_a^*
s_a ) .
\]

\noindent Hence $s_a s_b$ is a partial isometry which is nonzero
if and only if $s_a$ and $s_b$ are. Therefore we may define $s_{a
b} = s_a s_b$ and similarly define $s_\alpha$ for all $\alpha \in
\mathcal{L} (E , \pi )$. One checks that Definition \ref{lgdef}
(ii) holds for $\alpha \in \mathcal{L} ( E , \pi )$, Definition
\ref{lgdef} (iii) holds for $\alpha , \beta \in \mathcal{L}^n ( E
, \pi )$ for $n \ge 1$ and Definition \ref{lgdef} (iv) holds for
$A \in \mathcal{C}$ with finite and nonempty $L^n_A$ for $n \ge
1$. Then (cf.\ (\ref{mrels})) we have
\begin{equation} \label{Matrel}
s_\alpha^* s_\alpha s_\beta = p_{r ( \alpha )} s_\beta = s_\beta
p_{r ( r ( \alpha ) , \beta ) )} = s_\beta p_{r ( \alpha \beta )}
= s_\beta s_{\alpha \beta}^* s_{\alpha \beta} .
\end{equation}

\noindent To justify the requirement that $( E , \pi , \mathcal{C}
)$ is weakly left-resolving in Definitions \ref{lgdef}, consider
the following: Let $\{ p_A , s_a \}$ be a representation of $( E ,
\pi , \mathcal{C} )$ in which $p_A \ne 0$ for all $A \in
\mathcal{C}$. By Definition \ref{lgdef} (i) we have $(p_A - p_{A
\cap B})( p_B - p_{A \cap B} ) = 0$ for all $A, B \in
\mathcal{C}$. Suppose, for contradiction, that there is $\alpha
\in \mathcal{L} ( E , \pi )$ such that $r (A ,  \alpha ) \cap r (
B , \alpha ) \neq r ( A \cap B, \alpha)$. From Definition
\ref{lgdef} (iv) we have
\[
p_A - p_{A \cap B} \ge s_\alpha \left( p_{r ( A , \alpha )} -
p_{r(A \cap B,\alpha)} \right)s _\alpha^* \text{ and } p_B - p_{A
\cap B} \ge s_\alpha \left( p_{r ( B , \alpha )} - p_{r(A \cap B,
\alpha)}\right) s_\alpha^*
\]

\noindent so $(p_A - p_{A \cap B})( p_B - p_{A \cap B}) \neq 0$, a
contradiction. Thus a representation of $( E , \pi , \mathcal{C}
)$  will be degenerate if $( E , \pi , \mathcal{C} )$ is not
weakly left-resolving.

\noindent
Relation (iv) in Definition \ref{lgdef} can make sense
even if $A \in {\mathcal C}$ emits infinitely many edges in $E$: If there are only
finitely many different labels attached to the edges which $A$
emits then $L^1_A$ is finite. For directed graphs the analogue of
equation (\ref{sumcond}) holds when a vertex has finite valency; when this
is true at every vertex, the graph is called row-finite. With
this in mind, we make the following definition:

\begin{dfn}
Let $( E , \pi , \mathcal{C} )$ be a labelled space. We say
that $A \in {\mathcal C}$ is {\em singular} if $L^1_A$ is infinite.  If
no set $A \in {\mathcal C}$ is singular we say that $( E , \pi , \mathcal{C} )$
is {\em set-finite}.
\end{dfn}

\noindent If $( E , \pi , \mathcal{C} )$ is set-finite, then
$\displaystyle L_A^n$ is finite for all $A \in \mathcal{C}$ and
all $n \ge 1$. In the examples below, the resulting labelled space
will be set-finite whenever the original graph is row-finite.

\begin{examples} \label{E0ex}
\begin{itemize}

\item[(i)] Let $E$ be a directed graph with the trivial labelling
$\pi$. Then $\mathcal{E}^0$ consists of all
the finite subsets of $E^0$. If $E$ is row-finite then $( E , \pi , \mathcal{E}^0 )$
and $( E , \pi , \mathcal{E}^0_- )$ are set-finite. One may show that a representation of
$( E , \pi , \mathcal{E}^0 )$
is a Cuntz-Krieger $E$-family and conversely (see
\cite{bhrsz,bprsz} for instance).  If all sources in $E$
have finite valency, then the $*$-algebra generated by a
representation of $( E , \pi , \mathcal{E}^0_-
)$ contains a representation of $( E , \pi ,
\mathcal{E}^0 )$. If there is a source $v \in E^0$ with infinite valency
then there is no representative of $p_v$ in the $*$-algebra generated by a representation
of $( E , \pi , \mathcal{E}^0_- )$.

\item[(ii)] Under the identification of an ultragraph
$\mathcal{G}$ with a labelled graph $(E_{\mathcal G} ,
\pi_{\mathcal G} )$ we have $\mathcal{E}^0_{\mathcal{G}} =
\mathcal{G}^0$. Since $\mathcal{A} = \mathcal{G}^1$ a
representation of $( E_{\mathcal{G}} , \pi_{\mathcal G} ,
\mathcal{E}^0_{\mathcal{G}} )$ is a Cuntz-Krieger
$\mathcal{G}$-family (see \cite[Definition 2.7]{t}). If
$\mathcal{G}$ has sources which are singular then we get similar
behaviour to that described in (i) above.

\item[(iii)] In Examples \ref{lgex} (iii) we have $\mathcal{E}_i^0 ,
 ( \mathcal{E}_i^0 )_- = 2^{E_i^0}$ for $i=1,2$ whereas $\{ w \} \not\in ( \mathcal{E}_3^0 )_-$.
 A representation of $( E_2 , \pi_2 , ( \mathcal{E}_2^0)_- )$, is generated by partial isometries
$s_0 , s_1$ satisfying the relations in \cite[Proposition
8.3]{m} and \cite[\S 2]{c} for $\mathcal{O}_Y$, where $Y$ is the even shift. 

\item[(iv)] A covering $p : F \to E$ of directed graphs yields a
labelling $\pi_p : F^1 \to E^1$. We may identify $\mathcal{F}^0$
with the collection of inverse images of the finite subsets of
$E^0$. A representation of $( F , \pi_p , \mathcal{F}^0 )$ is a
Cuntz-Krieger $E$-family. If $F$ has  sources with infinite
valency, then we get similar behaviour to that described in (i)
above.

\item[(v)] An outsplitting $E_s (
\mathcal{P} )$ of $E$ gives rise to a labelling $\pi : E_s (
\mathcal{P} )^1 \to E^1$. If $\mathcal{P}$ is proper
then we may identify
$\mathcal{E_{{}_{\mathcal{S}}} ( \mathcal{P} )}^0$ with the
collection of finite subsets of $E^0$, and a representation of
$( E_s ( \mathcal{P} ) , \pi , \mathcal{E_{{}_{\mathcal{S}}} ( \mathcal{P} )}^0 )$
is a Cuntz-Krieger $E$-family. If $E$ has sources with infinite valency
then, we get similar behaviour to that described in (i) above, even when
the outsplitting is proper.

\item[(vi)] An arbitrary shift $\Lambda \subseteq \mathcal{A}^{\bf
Z}$ gives rise to a left-resolving labelled graph $(E_\Lambda ,
\pi_\Lambda )$ with no sources or sinks. If $\mathcal{A}$ is
finite then the generators of $\mathcal{O}_\Lambda$ form a
representation of $( E_\Lambda , \pi_\Lambda , (
\mathcal{E}_\Lambda^0)_- )$ (cf.\ \cite{cm,m}).

\item[(vii)] An arbitrary shift $\Lambda \subseteq
\mathcal{A}^{\bf Z}$ gives rise to a left-resolving labelled graph
$(E_{\Lambda^*} , \pi_{\Lambda^*} )$ with no sources or sinks. If
$\mathcal{A}$ is finite then the generators of
$\mathcal{O}_{\Lambda^*}$ form a representation of $(
E_{\Lambda^*} , \pi_{\Lambda^*} , ( \mathcal{E}_{\Lambda^*}^0)_-
)$ (cf.\ \cite{cm,m}).
\end{itemize}
\end{examples}

\noindent Examples \ref{E0ex} (i)-(v)  show that it is possible
for $\mathcal{E}^0$ and $\mathcal{E}^0_-$ to be different, but for
the $*$-algebras generated by representations of $( E , \pi ,
\mathcal{E}^0 )$ and $( E , \pi , \mathcal{E}^0_- )$ to be the
same.

Let $(E , \pi , \mathcal{C} )$ be a labelled space. Let
$\mathcal{C}^* = \mathcal{L} ( E , \pi ) \cup \mathcal{C}$ and
extend $r,s$ to $\mathcal{C}^*$ by $r ( A ) = A$, $s(A)=A$ for all
$A \in \mathcal{C}$. For $A \in \mathcal{C}$, put $s_A = p_A$, so
$s_\beta$ is defined for all $\beta \in \mathcal{C}^*$.



\begin{lem} \label{dense}
Let $(E, \pi , \mathcal{C})$ be a weakly left-resolving labelled
space and $\{ s_{a} , p_A \}$ a representation of $(E, \pi ,
\mathcal{C} )$.  Then any nonzero product of $s_{a}$, $p_A$ and
$s_{b}^*$ can be written as a finite combination of elements of
the form $s_\alpha p_A s_\beta^*$ for some $A \in {\mathcal C}$,
and $\alpha, \beta \in \mathcal{C}^*$ satisfying $A \subseteq r (
\alpha) \cap r ( \beta)  \neq \emptyset$.
\end{lem}

\begin{proof}
Since $s_{\alpha}p_A s_{\beta}^*  = s_\alpha p_{r ( \alpha ) \cap A \cap r
( \beta )} s_\beta^* $ it follows that $s_\alpha p_A s_\beta^*$ is
zero unless $A \cap r ( \alpha ) \cap r ( \beta ) \neq \emptyset$ and without loss
of generality we may assume that $A \subseteq r ( \alpha ) \cap r ( \beta )$.
For $\alpha, \beta, \gamma,\delta \in {\mathcal L}(E,\pi)$ and
$A,B \in {\mathcal C}$ we have
\begin{equation} \label{defprods}
\left(s_{\alpha} p_A s_{\beta}^*\right) \left(s_{\gamma} p_B
s_{\delta}^*\right) =
\begin{cases}
s_{\alpha\gamma'}p_{r(A,\gamma')\cap B } s_{\delta}^* & \mbox{ if } \gamma = \beta\gamma' \\
s_{\alpha}p_{A\cap   r(B,\beta')} s_{\delta\beta'}^* & \mbox{ if } \beta= \gamma\beta' \\
s_{\alpha}p_{  A \cap B } s_{\delta}^* & \mbox{ if } \beta = \gamma\\
0 & \mbox{ otherwise }
\end{cases}
\end{equation}

\noindent To see this, suppose $\gamma = \beta \gamma'$ then as $A
\subseteq r ( \beta ) \cap r ( \alpha )$

\begin{eqnarray*} s_\alpha p_A s_\beta^* s_\gamma p_B s_\delta^* &=& s_\alpha
p_A s_{\beta}^* s_\beta s_{\gamma'} p_B s_\delta^* = s_\alpha p_A
p_{r ( \beta)} s_{\gamma'} p_B s_\delta^*\\
 &=& s_\alpha p_A
s_{\gamma'} p_B s_\delta^* = s_{\alpha \gamma'} p_{r ( A , \gamma'
) \cap B} s_\delta^*
\end{eqnarray*}

\noindent A similar calculation gives the desired formulas in the
cases $\beta = \gamma \beta'$ and $\beta = \gamma$. If $\beta$ and
$\gamma$ have no common initial segment, then without loss of
generality, assume that $\beta \in \mathcal{L}^n (E,\pi)$ and
$\gamma \in \mathcal{L}^m (E,\pi )$ with $n > m$. Write $\beta=
\beta' \beta''$ where $\beta' \in \mathcal{L}^m (E,\pi)$, and then
by Definition \ref{lgdef}(iv) we have
$s_\beta^* s_\gamma = s_{\beta''}^* s_{\beta'}^* s_\gamma = 0$
since $\beta' \neq \gamma$ and so $s_\alpha p_A
s_\beta^* s_\gamma p_A s_\delta^* =0$.
By Definition \ref{lgdef} (i) and (ii) we may extend
(\ref{defprods}) to the case when $\alpha, \beta, \gamma , \delta
\in \mathcal{C}^*$.
\end{proof}

\begin{thm} \label{existence}
Let $( E, \pi , \mathcal{C} )$ be a weakly left-resolving labelled space.
There exists a
$C^*$-algebra $B$ generated by a universal representation of
$\{s_a, p_A \}$ of $( E , \pi  , \mathcal{C})$.  Furthermore the $s_a$'s are
nonzero and every $p_A$ with $A \neq \emptyset$ is nonzero.
\end{thm}

\begin{proof}
Let $S_{(E,\pi , \mathcal{C})} := \{ ( \alpha , A , \beta ): \alpha , \beta \in
\mathcal{C}^* , A \in {\mathcal C}, A \subseteq r ( \alpha )
\cap r( \beta ) \}$ and let $k_{( E , \pi , \mathcal{C} )}$ be the space of
functions of finite support on $S_{( E , \pi , \mathcal{C} )}$.  The set of
point masses $\{ e_{\tau} : \tau \in S_{(E,\pi , \mathcal{C} )} \}$ forms a
basis for $k_{( E , \pi , \mathcal{C} )}$.  Set $( \alpha , A , \beta )^* := (
\beta , A , \alpha )$; then thinking of $e_{( \alpha , A , \beta
)}$ as $s_{\alpha} p_A s_{\beta}^*$ and using (\ref{defprods}) we
can define a multiplication with respect to which $k_{(E,\pi , \mathcal{C} )}$ is
a $*$-algebra.

As a $*$-algebra $k_{( E , \pi , \mathcal{C} )}$ is generated by
$q_A := e_{(A,A,A)}$ for $A \in \mathcal{C}$ and $t_{a} := e_{( a
, r ( a ), r( a ) )}$ for $a \in \mathcal{L}^1 ( E , \pi )$. Our
definition of multiplication ensures that properties (ii) and
(iii) of Definition \ref{lgdef} hold; moreover $q_A q_B = q_{ A
\cap B }$. We mod out by the ideal $J$ generated by the elements
$q_{A \cup B} - q_A - q_B + q_{A \cap B}$ for $A, B \in
\mathcal{C}$, and $q_A - \sum_{a \in L_A^1 } s_{a }p_{r( A , a
)}s_{a }^*$ for $A \in \mathcal{C}$ with $L_A^1$ nonempty and
finite.  Then the images $r_A$ of $q_A$ and $u_{a}$ of $t_{a}$ in
$k_{( E, \pi , \mathcal{C} )}/J$ form a representation of $( E ,
\pi , \mathcal{C} )$ that generates $k_{(E , \pi , \mathcal{C}
)}/J$.  The triple $(k_{(E, \pi , \mathcal{C} )} / J , r_A ,
u_{a})$ has the required universal property, but is not a
$C^*$-algebra. Using a standard argument we can convert this
triple to a $C^*$-algebra $B$ satisfying the required properties
(see \cite[Theorem 2.1]{ahr} for instance).

Now for each $a \in \mathcal{L}^1 (E, \pi )$ and
$e \in \pi^{-1}(a)$, let $\mathcal{H}_{(a,e)}$ be an
infinite-dimensional Hilbert space. Also for each $v \in s (a)$ we
define ${\mathcal H}_{(a,v)} := \oplus_{\{ e : s(e)=v , \pi(e)=a
\}} {\mathcal H}_{(a,e)}$. If $v$ is a sink let ${\mathcal H}_v$
be an infinite-dimensional Hilbert space. For $A \in {\mathcal
C}$ we define ${\mathcal H}_A := \oplus_{b \in L^1_A} \oplus_{v
\in s( b ) \cap A} {\mathcal H}_{( b , v)}$ and then note that
each Hilbert space we have defined is a subspace of
\[
{\mathcal H} := \left( \oplus_{a \in \mathcal{L}^1 (E,\pi)}
\oplus_{v \in s(a)} {\mathcal H}_{(a,v)} \right) \oplus_{\{ v :
s^{-1} (v) = \emptyset \}} \mathcal{H}_v .
\]

%

\noindent For each $a \in {\mathcal L}^1(E,\pi)$, let $S_{a}$ be a
partial isometry with initial space ${\mathcal H}_{r(a)}$ and
final space $\oplus_{v \in s (a)} {\mathcal H}_{(a,v)} \subseteq
{\mathcal H}_{s ( a )}$. For $A \in {\mathcal C}$, define
$P_A$ to be the projection of ${\mathcal H}$ onto ${\mathcal
H}_A$, where this is interpreted as the zero projection when $A =
\emptyset$.

It is easy to verify that since $( E , \pi , \mathcal{C} )$ is
weakly left-resolving, the operators $\{S_{a},P_A\}$ form a
representation of $( E, \pi, \mathcal{C})$ in which $S_a , P_A$ are nonzero.
By the universal property there
exists a homomorphism $\pi_{S,P} : B \to C^*( \{ S_{a} , P_A \}
)$. Since the $S_{a}$'s and $P_A$'s are nonzero, it follows that
the $s_a$'s and $p_A$'s are also nonzero.
\end{proof}

\begin{dfn}
Let $( E , \pi , \mathcal{C} )$ be a weakly left-resolving
labelled space, then $C^* ( E , \pi , \mathcal{C} )$ is the
universal $C^*$-algebra generated by a representation of $( E , \pi , \mathcal{C} )$.
\end{dfn}

\noindent Let $( E , \pi , \mathcal{C} )$ be a weakly left-resolving
labelled space and $\{ s_a , p_A \}$ be the universal representation
of $( E , \pi , \mathcal{C} )$, then by Lemma \ref{dense}
\[
\text{span} \, \{ s_\alpha p_A s_\beta^* : \alpha , \beta \in \mathcal{L} ( E , \pi )
, A \in \mathcal{C} , A \subseteq r ( \alpha ) \cap r ( \beta ) \}
\]

\noindent is a dense $*$-subalgebra of $C^* ( E , \pi , \mathcal{C} )$.
 The following result may be proved along the same lines
as \cite[Lemma 3.2]{t}.

\begin{lem} \label{unital}
Let $\mathcal{A}$ be finite, $E$ have no sinks, and $( E , \pi , \mathcal{C} )$ be a weakly
left-resolving labelled space. Then $C^* ( E , \pi , \mathcal{C} )$ is unital.
\end{lem}

\begin{proof}
Observe that $\sum_{a \in \mathcal{A}} s_a s_a^* $ is a unit for $C^* ( E , \pi , \mathcal{C} )$.
\end{proof}


\begin{lem} \label{lgiso}
If $\phi : ( E, \pi ) \to (F, \pi')$ is  a labelled graph
isomorphism, then for all ${\mathcal C}$ which are accommodating
for $(E,\pi)$ we have $C^* (E, \pi , \mathcal{C} ) \cong C^* ( F ,
\pi' , \phi ( \mathcal{C} ) )$.
\end{lem}

\begin{proof}
The map $\phi$ induces a bijection between the generators of $C^*
(E , \pi , \mathcal{C} )$ and $C^* ( F , \pi' , \phi ( \mathcal{C} ) )$ and so by the universal
property there are homomorphisms from one $C^*$-algebra to the
other which are also inverses of each other.
\end{proof}

\section{Gauge Invariant Uniqueness Theorem} \label{GIUT}

\noindent Let $\{ s_a , p_A \}$ be the universal representation of
$( E, \pi , \mathcal{C} )$ which generates $C^* ( E, \pi , \mathcal{C} )$. For $z \in {\bf T}$,
$a \in \mathcal{L}^1 ( E, \pi )$ and $A \in \mathcal{C}$ let
\[
t_a := \gamma_z s_a = z s_a \text{ and } q_A := \gamma_z p_A = p_A
\]

\noindent then the family $\{ t_a , q_A \} \in C^* ( E, \pi , \mathcal{C} )$ is
also a representation of $( E, \pi , \mathcal{C} )$. By universality of $C^*( E, \pi , \mathcal{C} )$
and a routine $\epsilon/3$ argument we see that $\gamma$ extends
to  a strongly continuous action
$$
\gamma: {\bf T} \to \operatorname{Aut}\; C^* ( E, \pi , \mathcal{C} )
$$
which we call the {\em gauge action}.

\begin{prop}
\begin{itemize}
\item[i)] Let $E$ be a directed graph with the trivial labelling
$\pi$. Then $C^*(E,\pi, \mathcal{E}^0) \cong C^*(E)$.

\item[ii)] Let ${\mathcal G}$ be an ultragraph. Then
$C^*(E_{\mathcal{G}},\pi_{\mathcal G}, {\mathcal E}_{\mathcal
G}^0) \cong C^*({\mathcal G})$, where $(E_{\mathcal
G},\pi_{\mathcal{G}})$ is the labelled graph associated to
$\mathcal{G}$ .

\item[iii)] Let $p: F \to E$ be a covering map with induced
labelling $\pi_p: F^1 \to E^1$.  Then $C^*(F,\pi_p,{\mathcal F}^0)
\cong C^*(E)$.

\item[iv)] Let $E$ be a directed graph and let $E_s({\mathcal P})$
be an outsplitting. Let $\pi$ be the labelling of $E_s({\mathcal
P})$ induced by the outsplitting.  If ${\mathcal P}$ is a proper
partition then $C^*(E_s({\mathcal P}), \pi, {\mathcal
E}_s({\mathcal P})^0) \cong C^*(E)$.
\end{itemize}
\end{prop}

\begin{proof} In each case the left hand side contains a generating set for the
$C^*$-algebra on the right as shown in Examples \ref{E0ex}.  We apply the
appropriate gauge-invariant uniqueness theorem for the algebra on the right hand
side to obtain the isomorphism.
\end{proof}

\noindent To establish connections with the Matsumoto algebras we
need a version of the gauge-invariant uniqueness theorem for
labelled graph algebras.

\begin{lem} \label{keythought}
Let $( E , \pi , \mathcal{C} )$ be a weakly left-resolving labelled space, $\{ s_{a} ,
p_A \}$ a representation of $( E, \pi, \mathcal{C} )$, and $Y = \{ s_{\alpha_i}
p_{A_i} s_{\beta_i}^* : i = 1 , \ldots , N \}$ be a set of partial
isometries in $C^* ( E, \pi, \mathcal{C} )$ which is closed under multiplication
and taking adjoints. If $q$ is a minimal projection in $C^* (Y)$
then either
\begin{itemize}
\item[(i)] $q = s_{\alpha_i} p_{A_i} s_{\alpha_i}^*$ for some $1
\le i \le N$

\item[(ii)] $q = s_{\alpha_i} p_{A_i} s_{\alpha_i}^* - q'$ where
$q' = \sum_{l=1}^m s_{\alpha_{k(l)}} p_{A_{k(l)}}
s_{\alpha_{k(l)}}^*$ and $1 \le i \le N$; moreover there is a
nonzero $r = s_{\alpha_i \beta} p_{r ( A_i , \beta ) } s_{\alpha_i \beta}^* \in
C^* ( E, \pi, \mathcal{C} )$ such that $q' r =0$ and $q \ge r$.
\end{itemize}
\end{lem}

\begin{proof}


By equation (\ref{defprods}) any projection in $C^* (Y)$ may be
written as
\[
\sum_{j=1}^n s_{\alpha_{i(j)}} p_{A_{i(j)}} s_{\alpha_{i(j)}}^* -
\sum_{l=1}^m s_{\alpha_{k(l)}} p_{A_{k(l)}} s_{\alpha_{k(l)}}^*
\]

\noindent where the projections in each sum are mutually
orthogonal and for each $l$ there is a unique $j$ such that
$s_{\alpha_{i(j)}} p_{A_{i(j)}} s_{\alpha_{i(j)}}^* \ge
s_{\alpha_{k(l)}} p_{A_{k(l)}} s_{\alpha_{k(l)}}^*$.

If $q = \sum_{j=1}^n s_{\alpha_{i(j)}} p_{A_{i(j)}}
s_{\alpha_{i(j)}}^* - \sum_{l=1}^m s_{\alpha_{k(l)}} p_{A_{k(l)}}
s_{\alpha_{k(l)}}^*$ is a minimal projection in $C^* (Y)$ then we
must have $n=1$. If $m=0$ then $q = s_{\alpha_i} p_{A_i}
s_{\alpha_i}^*$ for some $1 \le i \le N$. If $m \neq 0$ then $q =
s_{\alpha_i} p_{A_i} s_{\alpha_i}^* - q'$ where $q' = \sum_{l=1}^m
s_{\alpha_{k(l)}} p_{A_{k(l)}} s_{\alpha_{k(l)}}^*$ and $1 \le k
\le N$. Since $q'$ is the sum of finitely many projections and $q
\neq 0$ it follows by repeated use of Definition \ref{lgdef} (iv)
that there is a nonzero $r = s_{\alpha_i \beta} p_{r ( A_i , \beta
)} s_{\alpha_i \beta}^*$ in $C^* (E, \pi, \mathcal{C} )$ such that
$r q' =0$ and $q \ge r$.
\end{proof}

\begin{thm} \label{giut}
Let $( E, \pi, \mathcal{C} )$ be a weakly left-resolving labelled space and let $\{ S_a
, P_A \}$ be a representation of $(E, \pi, \mathcal{C} )$ on Hilbert space.
Take $\pi_{S,P}$ to be the representation of $C^*( E, \pi, \mathcal{C} )$
satisfying
$\pi_{S,P}(s_a) = S_a$ and $\pi_{S,P}(p_A) = P_A$. Suppose that
each $P_A$ is non-zero whenever $A \neq \emptyset$, and that there
is a strongly continuous action $\beta$ of ${\bf T}$ on
$C^*(S_{\alpha},P_A)$ such that for all $z \in {\bf T}$, $\beta_z
\circ \pi_{S,P} = \pi_{S,P} \circ \gamma_z$. Then $\pi_{S,P}$ is
faithful.
\end{thm}

\begin{proof}
A straightforward argument along the lines of \cite[Lemma 2.2.3
]{pr1} shows that
\[
C^* (E , \pi, \mathcal{C} )^\gamma = \clsp \{ s_\alpha p_A s_\beta^* : \alpha ,
\beta \in \mathcal{L}^n (E,\pi) \text{ for some } n \text{ and } A
\subseteq r ( \alpha )
 \cap r ( \beta )  \}
\]

\noindent where $C^* ( E, \pi, \mathcal{C} )^\gamma$ is the fixed point algebra of
$C^* ( E, \pi, \mathcal{C} )$  under the gauge action $\gamma$. We claim that $C^*
(E , \pi , \mathcal{C} )^\gamma$ is AF. Let $Y$ be a finite subset of $C^* ( E , \pi , \mathcal{C}
)^\gamma$. Since $y \in Y$ may be approximated by a finite linear
combination of elements of the form $s_\alpha p_A s_\beta^*$ where
$\vert \alpha \vert = \vert \beta \vert$ we may assume that $Y =
\{ s_{\alpha_i} p_{A_i} s_{\beta_i}^* : \vert \alpha_i \vert =
\vert \beta_i \vert, i = 1 \ldots N \}$.

Let $M$ be the length of the longest word in $\{ \alpha_1 , \dots
, \alpha_N \}$. Let $W$ denote the collection of all words in
$\mathcal{L} (E,\pi)$ of length at most $M$ that can be formed
from composing subwords of $\alpha_1 , \dots , \alpha_N , \beta_1
, \ldots , \beta_N$. Let $\mathcal{B}$ be the collection all
finite intersections of $\{ A_i \}_{i=1}^n$ and $\{ r ( A_i ,
\gamma ) : 1 \le i \le N , \gamma \in W \}$.  By equation
(\ref{defprods}) a non-zero product of elements of $Y$ is of the
form $s_{\gamma} p_A s_{\delta}^*$ where $\gamma , \delta \in W$
and $A \in \mathcal{B}$. Since $W$ and $\mathcal{B}$ are finite it
follows that $Y' = \{ s_\gamma p_A s_\delta^* : \gamma , \delta
\in W , A \in \mathcal{B} \}$ is finite, closed under adjoints and
$C^* (Y) = C^*(Y')$. Hence we may assume that
$Y$ is closed under multiplication and taking adjoints. Thus $C^*
(Y) = \clsp (Y)$ is finite dimensional and so $C^* (E , \pi , \mathcal{C}
)^\gamma$ is AF by \cite[Theorem 2.2]{b}, establishing our claim.

To show that the canonical map $\pi_{S,P} : C^* (E , \pi,
\mathcal{C} ) \to C^* ( S_a , P_A )$ is injective on $C^* (E ,
\pi, \mathcal{C} )^\gamma$ we write $C^* (E , \pi, \mathcal{C}
)^\gamma$ as $\overline{\cup \; C^* ( Y_n )}$ where $\{ Y_n : n
\ge 1 \}$ is an increasing family of finite sets which are closed
under multiplication and taking adjoints. Suppose, for
contradiction, that $\pi_{S,P}$ is not faithful on $C^* ( Y_n )$
for some $n$. Then its kernel is an ideal and so must contain a
nonzero minimal projection $q$. If $Y_n = \{ s_{\alpha_i} p_{A_i}
s_{\beta_i}^* : i = 1 \ldots , N(n) \}$ then by Lemma
\ref{keythought} either $q = s_{\alpha_i} p_{A_i} s_{\alpha_i}^*$
for some $1 \le i \le N(n)$ or $q = s_{\alpha_i} p_{A_i}
s_{\alpha_i}^* - q'$ where $q' = \sum_{k=1}^m s_{\alpha_{i(k)}}
p_{A_{i(k)}} s_{\alpha_{i(k)}}^*$ and $1 \le i \le N(n)$. In the
first case $\pi_{S,P} ( s_{\alpha_i} p_{A_i} ) = S_{\alpha_i}
P_{A_i}$ is a partial isometry with initial projection $P_{A_i}$
and final projection $S_{\alpha_i} P_{A_i} S_{\alpha_i}^*$. But
$P_{A_i} = \pi_{S,P} ( p_{A_i} ) \neq 0$ by hypothesis and so
$\pi_{S,P} (q) = \pi_{S,P}( s_{\alpha_i} p_{A_i} s_{\alpha_i}^* )
= S_{\alpha_i} P_{A_i} S_{\alpha_i}^* \neq 0$ which is a
contradiction. In the second case by Lemma \ref{keythought} (ii)
there is $r = s_{\alpha_i \beta} p_{r ( A_i , \beta )} s_{\alpha_i
\beta}^*$ such that $q \ge r$ and $q' r =0$. We may apply the
above argument to show that $\pi_{S,P} (r) \neq 0$ and hence
$\pi_{S,P} (q) \ge \pi_{S,P} (r) \neq 0$ which is also a
contradiction. Hence $\pi_{S,P}$ is injective on $C^* ( Y_n )$ and
the result follows by arguments similar to those in \cite[Theorem
2.1]{bprsz}.
\end{proof}

\section{Applications} \label{App}

\subsection{Dual Labelled Graphs} \label{dualgraphs}

\noindent Let $E$ have no sinks and $(E, \pi )$ be a labelled
graph over alphabet $\mathcal{A}$. From this data we may form the
{\em dual labelled graph} $( \widehat{E} , \widehat{\pi} )$ over
alphabet $\widehat{\mathcal{A}} := \mathcal{L}^2 (E,\pi)$ as
follows: Let $\widehat{E}^0 = E^1$, $\widehat{E}^1 = E^2$ and the
maps $r', s': \widehat{E}^1 \to \widehat{E}^0$ be given by $r' (
ef ) = f \mbox{ and } s' ( ef ) = e$.  The labelling
$\widehat{\pi} : \widehat{E}^1 \to \widehat{\mathcal{A}}$ is
induced by the original labelling, so that $\widehat{\pi} ( ef ) =
\pi (e) \pi (f)$. For $ab \in \mathcal{L}^1 ( \widehat{E} ,
\widehat{\pi} ) = \mathcal{L}^2(E,\pi)$ we have
\[
r_{\widehat{\pi}} ( ab ) = \{ f : \widehat{\pi} ( ef ) = ab  \}
, \mbox{ and } s_{\widehat{\pi}} ( ab ) =
\{ e : \widehat{\pi} ( ef ) = ab  \}
\]

\noindent and for $B \in 2^{E^1}$
\[
r_{\widehat{\pi}} ( B , ab ) = \{ f : \pi (ef) = ab  ,  e \in B \} .
\]

\noindent These maps extend naturally to $\mathcal{L} ( \widehat{E} , \widehat{\pi} ) =
\cup_{n \geq 1} {\mathcal L}^n(\widehat{E},\widehat{\pi})$ where
for $n \geq 1$, ${\mathcal L}^n(\widehat{E},\widehat{\pi})$ is identified
with ${\mathcal L}^{n+1}(E,\pi)$.
Consider the following subsets of $2^{E^1}$
\begin{align*}
\widehat{\mathcal{E}} &= \{ \{ e \} : s(e) \text{ is a source} \} \cup
\{ r_{\widehat{\pi}} ( \alpha ) : \alpha \in \mathcal{L} ( \widehat{E} , \widehat{\pi} ) \}
\cup \{ s_{\widehat{\pi}} ( \alpha ) : \alpha \in \mathcal{L} ( \widehat{E}, \widehat{\pi} ) \} \\
\widehat{\mathcal{E}_-} &= \{ r_{\widehat{\pi}} ( \alpha ) : \alpha \in
\mathcal{L} ( \widehat{E} , \widehat{\pi} ) \} .
\end{align*}

\noindent Let $\widehat{\mathcal{E}}^0$ (resp.\
$\widehat{\mathcal{E}}^0_-$) be the smallest collection of subsets
of $2^{E^1}$ containing $\widehat{\mathcal{E}}$ (resp.\
$\widehat{\mathcal{E}}_-$) which is accommodating for $(
\widehat{E} , \widehat{\pi} )$. One checks easily that if $( E,
\pi, \mathcal{C} )$ is left-resolving, then $( \widehat{E} ,
\widehat{\pi}, \widehat{\mathcal{C}} )$ is weakly left-resolving
for $\mathcal{C} = \mathcal{E}^0 , \mathcal{E}^0_-$.

\noindent
For $B \in \widehat{\mathcal{E}}^0$ (resp.\ $B \in \widehat{\mathcal{E}}^0_-$) we set
\[
\widehat{L}^1_B = \{ ab \in \mathcal{L}^1 (
\widehat{E} , \widehat{\pi} ) : s_{\widehat{\pi}} (ab) \cap B \neq
\emptyset \} .
\]

\noindent  If $E$ has no sources and sinks, the shift
$\textsf{X}_{(\widehat{E},\widehat{\pi})}$ determined by the dual
labelled graph $( \widehat{E} , \widehat{\pi} )$ of $( E , \pi )$
is the second higher block shift $\textsf{X}_{(E,\pi)}^{[2]}$
formed from $\textsf{X}_{(E,\pi)}$ (cf.\ \cite[\S 1.4]{lm}).

\begin{rmks} \label{thepointis}
Suppose that $ab \in \mathcal{L}^2 (E , \pi )$ then $c \in
L^1_{r(ab)}$ if and only if $bc \in
\widehat{L}^1_{r_{\widehat{\pi}} (ab)}$; moreover $r(r(ab),c) = r
( s ( r_{\widehat{\pi}} (ab) ) , bc )$. Suppose that $A \in
\mathcal{E}^0$ (resp. $A \in {\mathcal E}^0_-$) then $a \in L^1_A$ and
$ab \in \mathcal{L}^2 (E, \pi)$ if and only if
$ab \in \widehat{L}^1_{s^{-1} (A)}$.
\end{rmks}

\begin{thm} Let $(E,\pi)$ be a set-finite,
left-resolving labelled graph with no sinks then $C^* ( E , \pi , \mathcal{E}^0 )
\cong C^* ( \widehat{E} , \widehat{\pi} , \widehat{\mathcal{E}}^0 )$, moreover
$C^* ( E , \pi , \mathcal{E}^0_- )
\cong C^* ( \widehat{E} , \widehat{\pi} , \widehat{\mathcal{E}}^0_- )$.
\end{thm}

\begin{proof}
Let $\{ s_a , p_A \}$ be a representation of $( E , \pi ,
\mathcal{E}^0 )$ and $\{ t_{ab} , q_B \}$ be a representation of
$( \widehat{E} , \widehat{\pi} , \widehat{\mathcal{E}}^0 )$. For
$ab \in \mathcal{L}^1 ( \widehat{E} , \widehat{\pi} )$ and $B \in
\widehat{\mathcal{E}}^0$ let $T_{ab} = s_{a} s_{b} s_{b}^*$ and
\[
Q_B := \sum_{ab \in \widehat{L}^1_B} s_{ab}
p_{r(s(B),ab)}s_{ab}^*.
\]
\noindent Since $(E , \pi , \mathcal{E}^0 )$ is set-finite
$( \widehat{E} , \widehat{\pi} , \widehat{\mathcal{E}}^0 )$ is
set-finite by Remarks \ref{thepointis} and so the above sum is finite. One
checks that $\{ T_{ab} ,  Q_B \}$
is a representation of $( \widehat{E} , \widehat{\pi} , \widehat{\mathcal{E}}^0 )$.

By the universal property there is a homomorphism $\pi_{T,Q} : C^*
( \widehat{E} , \widehat{\pi} , \widehat{\mathcal{E}}^0 ) \to C^* ( E , \pi , \mathcal{E}^0)$ with
$\pi_{T,Q} ( t_{ab} ) = T_{ab}$ and $\pi_{T,Q} ( q_B ) = Q_B$.
Since $\pi_{T,Q}$ intertwines the respective gauge actions and
$Q_B \neq 0$ it follows from Theorem \ref{giut} that $\pi_{T,Q}$
is faithful. We claim that $\pi_{T,Q}$ is surjective.
For $a \in \mathcal{L}^1 ( E , \pi )$ we have
\begin{align*}
s_a = s_a p_{r(a)} = s_a \sum_{b \in L^1_{r(a)}} s_b p_{r ( r(a) ,
b )} s_b^* &= \sum_{b \in L^1_{r(a)}} s_a s_b s_b^* s_b p_{r(ab)}
s_b^* \\
&= \sum_{b \in L^1_{r(a)}} s_a s_b s_b^* \sum_{c \in L^1_{r(ab)
}} s_{bc} p_{r ( r(ab) , c )} s_{bc}^* \\
&= \sum_{b \in L^1_{r(a)}} T_{ab} \sum_{bc \in
\widehat{L}^1_{r_{\widehat{\pi}} (ab)}} s_{bc} p_{r(s(
r_{\widehat{\pi}} (ab)),bc)}s_{bc}^* \text{ by Remarks
\ref{thepointis} }
\\
&= \sum_{b \in L^1_{r(a)}} T_{ab} Q_{r_{\widehat{\pi}} (ab)}
\end{align*}
%
and so $s_a \in C^* ( T_{ab} , Q_B  )$. For $A \in \mathcal{E}^0$,
 by Remarks \ref{thepointis} we have
\begin{eqnarray*}
p_A &=& \sum_{a \in L^1_A} s_a p_{r(A,a)}s_a^*
    = \sum_{a \in L^1_A}s_a \sum_{b \in L^1_{r(A,a)}} s_b
    p_{r(r(A,a),b)} s_b^*s_a^* \\
    &=& \sum_{ab \in \widehat{L}^1_{s^{-1}(A)}} s_{ab} p_{r(A,ab)} s_{ab}^*
    = Q_{s^{-1}(A)}
\end{eqnarray*}
%
and so $p_A \in C^* ( T_{ab}, Q_B )$ which establishes our claim.
The second isomorphism is proved along similar lines.
\end{proof}

\subsection{Matsumoto Algebras} \label{ma}

\begin{thm} \label{wegetsomecmalgs}
Let $\Lambda$ be a shift space over a finite alphabet
$\mathcal{A}$ which satisfies condition $(I)$ and has left-Krieger
cover $(E_\Lambda,\pi_\Lambda)$ then $\mathcal{O}_\Lambda \cong
C^* ( E_\Lambda , \pi_\Lambda , ( \mathcal{E}_\Lambda)^0_- )$.
Moreover, if $\Lambda$ satisfies condition $(I^*)$, and has
predecessor graph $(E_{\Lambda^*},\pi_{\Lambda^*})$ then
$\mathcal{O}_{\Lambda^*} \cong
C^* ( E_{\Lambda^*} , \pi_{\Lambda^*} , ( \mathcal{E}_{\Lambda^*})^0_- )$.
\end{thm}

\begin{proof}

By definition every $A \in ( \mathcal{E}_\Lambda )^0_-$ can be
written as a union of sets of the form $A_j = \cap_{i=1}^{m(j)} r
( \mu_i^j )$ for $j=1 , \ldots , n$. For $\mu \in \Lambda^*$ let
$q_{r ( \mu )} = t_\mu^* t_\mu$, then since the projections $\{
t_\mu^* t_\mu : \mu \in \Lambda^* \}$ are mutually commutative
(see \cite[p.686]{m2}) we may define $q_{r ( \mu ) \cap r ( \nu )}
= q_{r ( \mu )} q_{r ( \nu )}$, and hence define $q_{A_j}$ for $1
\le j \le n$. By the inclusion-exclusion principle one may further
define
\[
q_A = \sum_{j=1}^n q_{A_j} - \sum_{j \neq k} q_{A_j} q_{A_k} + \ldots + (-1)^{n+1} q_{A_1} \cdots q_{A_n} .
\]

\noindent Using calculations along the lines of those in \cite[\S 3]{m} one
checks that $\{ t_a ,  q_A \}$ is a
representation of $( E_\Lambda , \pi_\Lambda , ( \mathcal{E}_\Lambda )^0_- )$.

Let $\{ s_a , p_A \}$ be a representation of $( E_\Lambda ,
\pi_\Lambda , ( \mathcal{E}_\Lambda )^0_- )$. By the universal
property for $C^* ( E_\Lambda , \pi_\Lambda , (
\mathcal{E}_\Lambda )^0_- )$ there is a map $\psi_{t,q} : C^* (
E_\Lambda , \pi_\Lambda , ( \mathcal{E}_\Lambda )^0_- ) \to
\mathcal{O}_\Lambda$ such that $\psi_{t,q} ( s_a ) = t_a$ and
$\psi_{t,q} ( p_A ) = q_A$, in particular $\psi_{t,q}$ is surjective.

Since $\Lambda$ satisfies condition $(I)$ it follows by Theorem
\ref{matsumotohasgauge} that $\mathcal{O}_\Lambda$ carries a
strongly continuous action $\beta$ of $\bf T$. Since $\beta_z
\circ \psi_{t,q} = \psi_{t,q} \circ \gamma_z$ for all $z \in {\bf
T}$ and $\psi_{t,q} ( p_A ) = q_A \neq 0$ it follows from Theorem
\ref{giut} that $\psi_{t,q}$ is injective, which completes the
proof of the first statement.

The second statement is proved similarly.
\end{proof}

\begin{rmks}
\begin{itemize}
\item[(i)] In \cite[\S 5]{cm} a condition (*) is given under which
for shift spaces $\Lambda$ satisfying (*) conditions $(I)$ and
$(I^*)$ are equivalent and ${\mathcal O}_{\Lambda} \cong {\mathcal
O}_{\Lambda^*}$.
  This suggests that if $\Lambda$ satisfies (*) then $( E_\Lambda ,
\pi_\Lambda  )$ is labelled graph isomorphic to
$( E_{\Lambda^*} ,\pi_{\Lambda^*})$ and the isomorphism
of ${\mathcal O}_\Lambda$ and ${\mathcal O}_{\Lambda^*}$ can be deduced from Theorem \ref{lgiso}.
However \cite[Theorem 6.1]{cm} shows that, in general, ${\mathcal O}_\Lambda$ and ${\mathcal O}_{\Lambda^*}$
are not isomorphic. In particular,
$( E_\Lambda ,\pi_\Lambda )$ and
$( E_{\Lambda^*} ,\pi_{\Lambda^*})$
are not labelled graph isomorphic in general.
\item[(ii)]
The isomorphism  of $C^* ( E_\Lambda , \pi_\Lambda , ( \mathcal{E}_\Lambda)^0_- )$
 and $\mathcal{O}_\Lambda$ identifies
$C^* ( p_A : A \in ( \mathcal{E}_\Lambda )^0_- )$ with
$A_\Lambda \subset \mathcal{O}_\Lambda$.
Recall from \cite[Corollary 4.7]{m2} that $A_\Lambda \cong C ( \Omega_\Lambda )$,
hence we may think of the elements of $( \mathcal{E}_\Lambda )^0_-$ as indexing closed sets in $\Omega_\Lambda$.
\end{itemize}
\end{rmks}

\subsection{ Finiteness Conditions} \label{finitenessconditions}

\begin{dfn} A labelled graph $(E , \pi )$ is {\em label-finite} if
$\vert \pi^{-1} ( a ) \vert < \infty$ for all $a \in \mathcal{L}^1
(E,\pi)$.
\end{dfn}

\noindent If $(E , \pi )$ is label-finite then $\pi^{-1} ( \alpha
)$ is finite for all $\alpha \in {\mathcal L} ( E , \pi )$ and so
all sets in $\mathcal{E}^0$ are finite (and conversely). If $(E ,
\pi )$ is label-finite then $( \widehat{E} , \widehat{\pi} )$ is
label-finite. If $E$ is row-finite and $(E , \pi )$ is label-finite then
$( E , \pi , \mathcal{E}^0 )$ is set-finite.

The following result generalises \cite[Corollary 2.5]{bprsz} (see
also \cite[Remark 3.3 (i)]{bp}).

\begin{thm}  \label{soficcase}
Let $( E, \pi  )$ be a row-finite left-resolving labelled graph
which is label-finite and satisfies $\{v\} \in {\mathcal E}^0$ for
all $v \in E^0$. Then $C^*( E, \pi , \mathcal{E}^0 ) \cong
C^*(E)$; moreover if $\{ v \} \in \mathcal{E}^0_-$ for all $v \in
E^0$ then $C^*( E, \pi , \mathcal{E}^0_- ) \cong C^*(E)$.
\end{thm}

\begin{proof}
Let $\{ s_e, p_v \}$ be the canonical Cuntz-Krieger $E$-family and
$\{ t_a , q_A \}$ be the canonical generators of $C^* ( E , \pi , \mathcal{E}^0
)$. For $a \in \mathcal{L}^1 ( E, \pi )$ and $A \in {\mathcal E}^0$
let
$$
T_a = \sum_{e \in E^1 : \pi ( e ) = a } s_{e} , \text{ and } Q_A =
\sum_{v \in A} p_v .
$$

\noindent The above sums make sense since $( E, \pi )$ is
label-finite. Since $E$ is row-finite one may easily check that
these operators define a representation of $( E , \pi ,
\mathcal{E}^0 )$. 
By the universal property of $C^* ( E, \pi , \mathcal{E}^0 )$ there is a
homomorphism $\psi_{T,Q} : C^* (E , \pi , \mathcal{E}^0 ) \to C^* (E)$ given by
$\psi_{T,Q} ( t_a ) = T_a$ and $\psi_{T,Q} ( q_A ) = Q_A$ for all $a
\in \mathcal{L}^1 ( E , \pi )$ and $A \in \mathcal{E}^0$.

Since $\{ v \} \in \mathcal{E}^0$ for all $v \in E^0$, we have
$p_v = Q_v \in C^*(  T_a , Q_A )$ for all $v \in E^0$. Since
our labelled graph is left-resolving we have $s_e = T_{\pi (e)}
Q_{r(e)} \in C^* (  T_a , Q_A )$ for all $e \in E^1$, and so
$\psi_{T,Q}$ is surjective. The canonical gauge actions  on
$C^*(E)$ and $C^* ( E, \pi , \mathcal{E}^0 )$ satisfy the required properties and
$\psi_{T,Q} ( q_A ) = Q_A \neq 0$ for all $A \in \mathcal{E}^0$, so
$\psi_{T,Q}$ is an isomorphism by Theorem \ref{giut}.

The proof of the second isomorphism is essentially the same.
\end{proof}

\begin{cor}
Let $\mathcal{G} = ( G^0 , \mathcal{G}^1 , r , s)$ be a row-finite
ultragraph then $C^* ( \mathcal{G} ) \cong C^* ( E_{\mathcal{G}}
)$ where $E_\mathcal{G}$ is the underlying directed graph of
$\mathcal{G}$.
\end{cor}

\begin{proof}
From Examples \ref{lgex} (ii) a row-finite ultragraph
$\mathcal{G}$ may be realised as a row-finite left-resolving
labelled graph $( E_{\mathcal{G}} , \pi_{\mathcal{G}} )$. As
$E_\mathcal{G}$ is row-finite it follows that $( E_{\mathcal{G}} ,
\pi_{\mathcal{G}} )$ is label-finite. Since the source map is
single-valued it follows that $v \in \mathcal{E}_{\mathcal{G}}^0$
for all $v \in G^0 = E_{\mathcal{G}}^0$ and hence the result
follows from Theorem \ref{soficcase}.
\end{proof}

\noindent The following result was first observed in \cite[Theorem
3.5]{c} (see also \cite[Theorem 4.4.4]{s}).

\begin{cor}
Let $\Lambda$ be a sofic shift over a finite alphabet then
\[
\mathcal{O}_\Lambda \cong C^* ( E_\Lambda )
\]

\noindent where $( E_\Lambda , \pi_\Lambda )$ is the left-Krieger cover of $\Lambda$.
\end{cor}

\begin{proof}
As ${E}_\Lambda^0$ is finite and each $v \in E_\Lambda^0$ has a
different past there are $\alpha_v \in \mathcal{L} (E_\Lambda ,
\pi_\Lambda )$ with $r_{\pi_\Lambda} ( \alpha_v ) = \{ v \}$.
Hence $\{ v \} \in ( \mathcal{E}_\Lambda^0 )_-$ for all $v \in
E_\Lambda^0$. The result follows by Theorem \ref{soficcase}.
\end{proof}

\noindent From \cite[Theorem 3.3.18]{lm} any two minimal
left-resolving representations $(E , \pi )$, $(F , \pi' )$ of an
irreducible sofic shift  are labelled graph isomorphic and so $C^*
( E, \pi , \mathcal{E}^0_- ) \cong C^* ( F, \pi' , \mathcal{F}^0_-
)$ by Lemma \ref{lgiso}. Moreover, one may use the minimality of
the representation to show that the underlying graph $E$ is
irreducible (cf.\ \cite[Lemma 3.3.10]{lm}). Hence we have:

\begin{cor} \label{soficsimple}
Let $( E , \pi )$ be a minimal left-resolving presentation of an
irreducible sofic shift over a finite alphabet, then $C^* ( E ,
\pi , \mathcal{E}^0_- ) \cong C^* ( E , \pi , \mathcal{E}^0 )$ is simple.
\end{cor}

\begin{rmk} Recall that the graph $( E_2 , \pi_2 )$ in Examples \ref{lgex} (ii)
is the left-Krieger cover of the even shift $Y$. Although $Y$ is
irreducible, $( E_2 , \pi_2 )$ is not a minimal left-resolving
presentation of $Y$ and $\mathcal{O}_Y \cong C^* ( E_2 )$ is not
simple. However the graph $( E_1 , \pi_1 )$ Examples \ref{lgex}
(ii) is a minimal left-resolving cover of $Y$ and so
\[
C^* ( E_1 , \pi_1 , \mathcal{E}^0_- ) \cong
C^* ( E_1 , \pi_1 , \mathcal{E}^0 ) \cong C^* ( E_1 )
\]

\noindent is simple. Similarly $C^* ( E_Z , \pi_Z ,
\mathcal{E}^0_- ) \cong C^* ( E_Z )$ is simple where $Z$ is the
irreducible shift introduced in Examples \ref{lgex} (vi).

Thus, if one wishes to associate a simple $C^*$-algebra to an
irreducible sofic shift $\Lambda$, then one should use the minimal
left-resolving presentation of $\Lambda$. (cf.\ \cite{c,c1}).
\end{rmk}

\noindent For a general shift space $\Lambda$, either $( E_\Lambda
, \pi_\Lambda )$ will not be row-finite or there will be $v \in
E_\Lambda^0$ with $v \not\in ( \mathcal{E}_\Lambda )^0_-$. This
indicates that the $C^*$-algebras corresponding to presentations
of such shift spaces will not be Morita equivalent to  graph
algebras. The shift associated to a certain Shannon graph (see
\cite[Theorem 7.7]{m5}) provides such an example.


\begin{thebibliography}{20}


\bibitem{bhrsz} T.~Bates, J-H.~Hong, I.~Raeburn and W.~Szyma\'{n}ski.
\newblock {\em The ideal structure of the $C^*$--algebras of infinite
graphs.}
\newblock Illinois J. Math {\bf 46} (2002), 1159--1176.

\bibitem{bprsz} T.~ Bates, D.~ Pask, I.~ Raeburn and W.~Szyma\'{n}ski.
\newblock {\em The $C^*$-algebras of row-finite graphs.}
\newblock New York J.\ Math. {\bf 6} (2000), 307--324.

\bibitem{bp} T.\ Bates and D.\ Pask. {\em Flow equivalence of graph
algebras}.  Ergod.\ Th.\ \& Dynam.\ Sys., {\bf 24} (2004),
367--382.

\bibitem{bp3} T.\ Bates and D.\ Pask. {\em Flow equivalence of labelled
graph algebras}.  In preparation.

\bibitem{b} O.~Bratteli. {\it Inductive limits of finite dimensional
${C}^*$-algebras}. Trans.\ Amer.\ Math.\ Soc., {\bf 171} (1972),
195--234.

\bibitem{c} T.\ Carlsen. {\em On $C^*$-algebras Associated with Sofic
Shifts}. J.\ Operator Theory {\bf 49} (2003), 203--212.

\bibitem{c1} T.\ Carlsen. {\em Symbolic dynamics, partial dynamical systems,
boolean algebras and $C^*$-algebras generated by partial
isometries}, preprint Univ.\ Oslo, (2004).

\bibitem{cm} T.\ Carlsen and K.\ Matsumoto. {\em Some remarks on
the $C^*$-algebras associated with subshifts}. Math.\ Scand.\ {\bf
95} (2004), 145--160.

\bibitem{dpr} K.\ Deicke,  D.~ Pask and I.~ Raeburn,
{\em  Coverings of directed graphs and crossed products of
$C^*$-algebras by coactions of homogeneous spaces}. Internat.\ J.\
Math., {\bf 14} (2003), 773-789.

\bibitem{ahr} A. an Huef and I. Raeburn. {\em The ideal structure of Cuntz-Krieger
 algebras}.  Ergod. Th. \& Dynam. Sys. {\bf 17} (1997), 611--624.

\bibitem{Kr} W.\ Krieger. {\em Sofic Systems I}.  Israel J.\ Math.,
{\bf 48}, (1984), 305-330.

\bibitem{kpr} A. ~Kumjian, D.~Pask and I.~Raeburn. {\em Cuntz-Krieger algebras of
directed graphs}. Pacific J. Math. {\bf 184} (1998), 161--174.

\bibitem{lm} D. Lind and B. Marcus. {\em An Introduction to Symbolic Dynamics and Coding}, CUP, 1995.

\bibitem{m} K.\ Matsumoto. {\em On $C^*$-algebras associated with
subshifts}. Internat.\ J.\ Math.\ {\bf 8}, (1997), 357-374.

\bibitem{m3} K.\ Matsumoto. {\em $K$-theory for  $C^*$-algebras associated with
subshifts}. Math.\ Scand.\ {\bf 82}, (1998), 237-255.

\bibitem{m1} K.\ Matsumoto. {\em Relations among generators of $C^*$-algebras associated with subshifts},
Internat.\ J.\ Math.\ {\bf 10} (1999), 385-405.

\bibitem{m2} K.\ Matsumoto. {\em Dimension groups for subshifts and simplicity
of the associated $C^*$-algebras}. J.\ Math.\ Soc.\ Japan {\bf 51}
(1999), 679-697.

\bibitem{m4} K.\ Matsumoto. {\em On automorphisms of  $C^*$-algebras associated with
subshifts}. J.\ Operator\ Theory {\bf 44}, (2000), 91-112.

\bibitem{m45} K.\ Matsumoto. {\em Bowen-Franks groups for
subshifts and Ext-groups for $C^*$-algebras}. K-Theory {\bf 23},
(2001), 67-104.

\bibitem{m5} K.\ Matsumoto. {\em Stabilized $C^*$-algebras associated with presentations
of subshifts}. Documenta\ Math.\ {\bf 7}, (2002), 1-30.

\bibitem{m6} K.\ Matsumoto. {\em Stabilized $C^*$-algebras constructed from symbolic
dynamical systems}. Ergodic\ Theory\ Dynam.\ Sys.\ {\bf 20},
(2000), 821-841.


\bibitem{pr1} D.~Pask and I.~Raeburn, \emph{On the {K}-theory of {C}untz-{K}rieger
algebras}.  Publ. RIMS Kyoto Univ. \textbf{32} (1996), 415--443.



\bibitem{rs} I.\ Raeburn and W.\ Szyma\'{n}ski. {\em Cuntz-Krieger algebras of
infinite graphs and matrices}. Trans.\ AMS, {\bf 356}, (2004)
39--59.

\bibitem{s} J.\ Samuel. {\em Graphs, Sofic Shifts and C*-algebras I}.
 Preprint Univ.\ Victoria, (1998).

\bibitem{Ta} R.G. Taylor. {\em Models of computation and formal
  languages}, Oxford University Press, 1998.

\bibitem{t} M. Tomforde. {\em A unified approach to Exel-Laca algebras and $C^*$-algebras associated
to graphs}. J. Operator Theory {\bf 50} (2003), 345--368.

\bibitem{t2} M. Tomforde. {\em Simplicity of ultragraph algebras}. Indiana Univ. Math. J. {\bf 52} (2003), 901--926.


\end{thebibliography}
\end{document}